\input amstex
\input amsppt.sty
\magnification=\magstep1
\parindent=1em
\baselineskip 15 pt
\hsize=12.3 cm
\vsize=18.5 cm
\nologo
\pageno=1
\TagsOnRight
\topmatter
\def\Z{\Bbb Z}

\def\Q{\Bbb Q}

\def\bg{\bigg}
\def\({\bg(}
\def\[{\bg[}
\def\){\bg)}
\def\]{\bg]}
\def\t{\text}
\def\f{\frac}
\def\colon{{:}\;}
\def\endrem{}

\def\Tor{\roman{Tor}}
\def\em{\emptyset}
\def\se {\subseteq}

\def\sm{\setminus}

\def\ls{\leqslant}
\def\gs{\geqslant}

\def\la{\lambda}

\def\colon{{:}\;}

\def\Proof{\noindent{\it Proof}}
\def\Remark{\medskip\noindent{\it  Remark}}
\def\Ack{\medskip\noindent {\bf Acknowledgments}}
\topmatter \hbox{Israel J. Math. 154(2006), no.\,1, 21--28.}
\medskip
\title Restricted sumsets and a conjecture of Lev\endtitle
\leftheadtext{Hao Pan and Zhi-Wei Sun}
\author Hao Pan and Zhi-Wei Sun*\endauthor
\affil Department of Mathematics, Nanjing University
\\ Nanjing 210093, People's Republic of China
\\ {\tt haopan79\@yahoo.com.cn}\ \ \quad\tt{zwsun\@nju.edu.cn}
\endaffil

\abstract Let $A,B,S$ be finite subsets of an abelian
group $G$. Suppose that the restricted sumset
$$C=\{a+b\colon a\in A,\ b\in B,\ \t{and}\ a-b\not\in S\}$$
is nonempty and some $c\in C$ can be written as $a+b$ with $a\in
A$ and $b\in B$ in at most $m$ ways. We
show that if $G$ is torsion-free or
elementary abelian then $|C|\gs|A|+|B|-|S|-m$.
We also prove that $|C|\gs|A|+|B|-2|S|-m$ if
the torsion subgroup of $G$ is cyclic. In the case $S=\{0\}$ this
provides an advance on a conjecture of Lev.
\endabstract
\thanks  \ \ 2000 {\it Mathematics Subject Classification}.
Primary 11B75; Secondary 05A05, 11P99, 20D60.
\newline\indent
*This author is responsible for communications, and supported by
the National Science Fund for Distinguished Young Scholars (no.
10425103) and the Key Program of NSF (no. 10331020) in China.
\endthanks
\endtopmatter

\document

\heading{1. Introduction}\endheading

Let $A$ and $B$ be finite nonempty subsets of an (additively written)
abelian group $G$. The sumset of $A$ and $B$ is defined by
$$A+B=\{a+b\colon a\in A\text{ and }b\in B\}.$$
The Cauchy-Davenport theorem (cf. [N, pp.\,43-48]),
a basic result in additive combinatorial number theory,
states that
$$|A+B|\gs\min\{p,\,|A|+|B|-1\}$$
if $G=\Z/p\Z$ with $p$ prime. Another theorem due to
Kemperman and Scherk (cf. [Sc], [Ke] and [L2]) asserts that
$$|A+B|\gs|A|+|B|-\min_{c\in A+B}\nu_{A,B}(c),\tag1.1$$
where
$$\nu_{A,B}(c)=|\{(a,b)\in A\times B\colon  a+b=c\}|;\tag1.2$$
in particular, we have $|A+B|\gs|A|+|B|-1$ if
some $c\in A+B$ can be uniquely written as $a+b$ with $a\in A$ and $b\in B$.

Now we define the restricted sumset
$$A\dotplus B=\{a+b\colon a\in A,\ b\in B,\ \t{and}\ a\not=b\}.\tag1.3$$
In 1964, Erd\H os and Heilbronn [EH] conjectured that
if $G=\Z/p\Z$ with $p$ prime then
$$|A\dotplus A|\gs\min\{p,\,2|A|-3\}.$$
This is much more difficult than the Cauchy-Davenport theorem
concerning unrestricted sumsets. It had been open for thirty years
until Dias da Silva and Hamidoune [DH] confirmed it in 1994
using representations of symmetric groups.
Later Alon, Nathanson and Ruzsa [ANR1, ANR2]
developed a powerful polynomial method to give a simpler proof of
the Erd\H os-Heilbronn conjecture (see also [A2]). They
showed that if $G=\Z/p\Z$ with $p$ prime then
$$|A\dotplus B|\gs \min\{p,\,|A|+|B|-2-\delta\},$$
where $\delta$ is $1$ or $0$ according to whether $|A|=|B|$ or not.
The reader may consult [HS], [K1], [K2], [L1], [LS] and [SY] for various
extensions of the Erd\H os-Heilbronn conjecture.

Motivated by the Kemperman-Scherk theorem
and the Erd\H os-Heilbronn conjecture,
Lev [L2] proposed the following interesting conjecture.

\proclaim{Conjecture 1.1 {\rm (Lev)}} Let $G$ be an
abelian group, and let $A$ and $B$ be finite nonempty subsets of
$G$. Then we have
$$|A\dotplus B|\gs|A|+|B|-2-\min_{c\in A+B}\nu_{A,B}(c).\tag1.4$$
\endproclaim

This conjecture is known to be true for torsion-free abelian groups
and elementary abelian $2$-groups. It also holds when $|G|$ is prime, or
$G$ is cyclic and $|G|\ls25$. (Cf. [L2].)

 Now we state our main results.

\proclaim{Theorem 1.1} Let $A$ and $B$ be finite nonempty subsets of a field $F$.
Let $P(x,y)\in F[x,y]$ and
$$C=\{a+b\colon  a\in A,\ b\in B, \ \t{and}\ P(a,b)\not=0\}.\tag1.5$$
If $C$ is nonempty, then
$$|C|\gs|A|+|B|-\deg P-\min_{c\in C}\nu_{A,B}(c).\tag1.6$$
\endproclaim

\Remark\ 1.1. When $P(x,y)=1$, (1.6) becomes (1.1).
\endrem
\medskip

Notice the difference between the minima in (1.4) and (1.6): as $C\se A+B$ we have
$\min_{c\in A+B}\nu_{A,B}(c)\ls\min_{c\in C}\nu_{A,B}(c)$.

\proclaim{Theorem 1.2}
Let $A$ and $B$ be finite nonempty subsets
of an abelian group $G$ whose torsion subgroup
$$\Tor(G)=\{g\in G\colon  g\ \t{has a finite order}\}$$
is cyclic.
For $i=1,\ldots,l$ let $m_i$ and $n_i$ be
nonnegative integers and let $d_i\in G$.
Suppose that
$$C=\{a+b\colon  a\in A,\ b\in B,\ \t{and}\ m_ia-n_ib\not=d_i
\ \t{for all}\ i=1,\ldots,l\}\tag1.7$$
is nonempty. Then
$$|C|\gs|A|+|B|-\sum_{i=1}^l(m_i+n_i)-\min_{c\in C}\nu_{A,B}(c).\tag1.8$$
\endproclaim

\Remark\ 1.2. When $A$ and $B$ are finite subsets of $\Z$,
the restricted sumset in (1.7) was first
studied by Sun [Su1].
\endrem
\medskip

From Theorems 1.1 and 1.2 we deduce the following result
on difference-restricted sumsets.

\proclaim{Theorem 1.3} Let $G$ be an abelian group, and
let $A,B,S$ be finite nonempty subsets of $G$ with
$$C=\{a+b\colon  a\in A,\ b\in B,\ \t{and}\ a-b\not\in S\}\not=\em.\tag1.9$$

{\rm (i)} If $G$ is torsion-free or elementary abelian, then
$$|C|\gs|A|+|B|-|S|-\min_{c\in C}\nu_{A,B}(c).\tag1.10$$

{\rm (ii)} If $\Tor(G)$ is cyclic, then
$$|C|\gs|A|+|B|-2|S|-\min_{c\in C}\nu_{A,B}(c).\tag1.11$$
\endproclaim

\Proof. Without loss of generality we can assume that $G$
is generated by the finite set $A\cup B\cup S$.

 If $G\cong \Z^n$, then we can simply view $G$ as
the ring of algebraic integers in an algebraic number field
$K$ with $[K:\Q]=n$. If $G\cong(\Z/p\Z)^n$ where $p$ is a prime,
then $G$ is isomorphic to the additive group of the finite field
with $p^n$ elements. Thus part (i) follows from Theorem 1.1 in the case
$P(x,y)=\prod_{s\in S}(x-y-s)$.

Let $d_1,\ldots,d_l$ be all the distinct elements of $S$.
Applying Theorem 1.2 with $m_i=n_i=1$ for all $i=1,\ldots,l$
we immediately get the second part. \qed

\Remark\ 1.3. It is interesting to compare Theorem 1.3 in the case $S=\{0\}$ with
Conjecture 1.1.
\endrem
\medskip

Concerning the set $C$ given by (1.9), there are some known results of different types.
When $A,B,S$ are finite nonempty subsets of
a field whose characteristic is an odd prime $p$,
the authors [PS] proved that
$|C|\gs\min\{p,\,|A|+|B|-|S|-q-1\}$,
where $q$ is the largest power of $p$ not exceeding $|S|$.
By modifying K\'arolyi's proof of [K1, Theorem 3], we can
show that if $q>1$ is a power of a prime $p$, and $A,B,S$ are subsets of $\Z/q\Z$ with
$\min\{|A|,|B|\}>|S|$, then $|C|\gs\min\{p,\,|A|+|B|-2|S|-1\}$.

We will give a key lemma in the next section and
prove Theorems 1.1 and 1.2 in Section 3. Our proofs use a version
of the polynomial method.

\heading{2. Some preparations}\endheading

Our basic tool is as follows.

\proclaim{Combinatorial Nullstellensatz {\rm ([A1, Theorem 1.1])}} Let
$A_1,\ldots,A_n$ be finite nonempty subsets of a field $F$, and set
$g_i(x)=\prod_{a\in A_i}(x-a)$ for $i=1,\ldots,n$.
Then $f(x_1,\ldots,x_n)\in F[x_1,\ldots,x_n]$
vanishes over the Cartesian product $A_1\times\cdots\times A_n$ if and only if
it can be written in the form
$$f(x_1,\ldots,x_n)=\sum_{i=1}^n g_i(x_i)h_i(x_1,\ldots,x_n)
$$
where $h_i(x_1,\ldots,x_n)\in F[x_1,\ldots,x_n]$ and
$\deg h_i\ls\deg f-\deg g_i$.
\endproclaim

With help of the Combinatorial Nullstellensatz, we provide
a lemma for our purposes.

\proclaim{Lemma 2.1}
Let $A$ and $B$ be finite nonempty subsets of a field $F$, and write
$$\nu_i=|\{(a,b)\in A\times B\colon a+\lambda_ib=\mu_i\}|\tag2.1$$
for $i=1,\ldots,k$ where $\lambda_i\in F\sm\{0\}$ and $\mu_i\in F$.
Let $P(x,y)\in F[x,y]$. Suppose that
for any $i=1,\ldots,k$ there are $a\in A$ and $b\in B$ with $P(a,b)\not=0$ and
$a+\la_i b=\mu_i$, and that
for each $(a,b)\in A\times B$ with $P(a,b)\not=0$ there is a unique $i\in\{1,\ldots,k\}$
with $a+\la_i b=\mu_i$.
Then we have
$$k+\min\{\nu_1,\ldots,\nu_k\}\gs|A|+|B|-\deg P.\tag2.2$$
\endproclaim
\Proof. Clearly
$$f(x,y):=P(x,y)\prod_{j=1}^k(x+\lambda_jy-\mu_j)$$
vanishes over $A\times B$. Set $g_A(x)=\prod_{a\in A}(x-a)$ and
$g_B(y)=\prod_{b\in B}(y-b)$. By the Combinatorial Nullstellensatz,
there are $h_A(x,y),h_B(x,y)\in F[x,y]$ such that
$$f(x,y)=g_A(x)h_A(x,y)+g_B(y)h_B(x,y)$$
and
$$\max\{\deg g_A+\deg h_A,\, \deg g_B+\deg h_B\}\ls \deg f.$$

Fix $1\ls i\ls k$. Write $h_B(x,y)=\sum_{s,t\gs0}c_{st}x^sy^t$ where $c_{st}\in F$.
Then
$$h_B(x,y)=\sum_{s,t\gs0}c_{st}((x+\la_i y-\mu_i)+\mu_i-\la_iy)^sy^t
=(x+\la_iy-\mu_i)q(x,y)+r(y),$$
where $q(x,y)\in F[x,y]$, and $r(y)=h_B(\mu_i-\la_iy,y)$ has degree not
greater than $\deg h_B$.

Now assume that
$k+\nu_i<|A|+|B|-\deg P$. We want to deduce a contradiction.
Set
$$A_0=\{a\in A\colon (\mu_i-a)/\lambda_i\not\in B\}.$$
Obviously $|A_0|=|A|-\nu_i$ and
$g_B((\mu_i-a)/\lambda_i)\not=0$ for any $a\in A_0$.
If $a\in A_0$, then
$$g_B\(\f{\mu_i-a}{\la_i}\)h_B\(a,\f{\mu_i-a}{\la_i}\)
=f\(a,\f{\mu_i-a}{\la_i}\)-g_A(a)h_A\(a,\f{\mu_i-a}{\la_i}\)=0$$
and hence
$$r\(\f{\mu_i-a}{\la_i}\)=h_B\(a,\f{\mu_i-a}{\la_i}\)=0.$$
Since $\deg r\ls\deg f-\deg g_B<|A|-\nu_i=|A_0|$, we must have
$r(y)=0$, i.e., $h_B(x,y)$ is divisible by $x+\lambda_iy-\mu_i$.
Recall that there are $a_0\in A$ and $b_0\in B$ such that
$P(a_0,b_0)\not=0$ and $a_0+\lambda_ib_0=\mu_i$. Since $h_B(a_0,b_0)=0$,
the polynomial
$P(a_0,y)\prod_{j=1}^k(a_0+\lambda_jy-\mu_j)=f(a_0,y)=g_B(y)h_B(a_0,y)$
is divisible by $(y-b_0)^2$.
As $a_0+\lambda_jb_0\not=\mu_j$ for any $j\not=i$,
we must have $y-b_0\mid P(a_0,y)$, which contradicts the fact that
$P(a_0,b_0)\not=0$. \qed

\heading{3. Proofs of Theorems 1.1--1.2}\endheading

\noindent{\it Proof of Theorem 1.1}.
Let $\mu_1,\ldots,\mu_k$ be all the distinct elements of $C$.
Applying Lemma 2.1 with $\la_1=\cdots=\la_k=1$, we find that
$$|C|+\min_{c\in C}\nu_{A,B}(c)\gs|A|+|B|-\deg P$$
which is equivalent to (1.6). \qed

\medskip
\noindent{\it Proof of Theorem 1.2}.
Without loss of generality, we can assume that $G$ is finitely generated,
and furthermore that $G$ is a subgroup of the multiplicative group
of the field of complex numbers (see the proof of Theorem 1.1 of [Su2]);
thus, $C$ is the set
$$\{ab\colon  a\in A,\, b\in B,
\ \t{and}\ a^{m_i}b^{-n_i}\not=d_i\ \t{for all}\ i=1,\ldots,l\}.$$

Let $-\la_1,\ldots,-\la_k$ be all the distinct elements of $C$, and
set
$$P(x,y)=\prod_{i=1}^l(x^{m_i}y^{n_i}-d_i).$$
Then, for each $j\in\{1,\ldots,k\}$, there are $a\in A$ and $b\in B$ such that
$a+\la_jb^{-1}=0$ and $P(a,b^{-1})\not=0$.
If $a\in A$, $b\in B$ and $P(a,b^{-1})\not=0$,
then there is a unique $j\in\{1,\ldots,k\}$
such that $\la_j=-ab$ (i.e., $a+\la_jb^{-1}=0$).
Applying Lemma 2.1 to the sets $A$ and $B^{-1}=\{b^{-1}\colon b\in B\}$
with $\mu_1=\cdots=\mu_k=0$,
we obtain that
$$k+\min_{1\ls j\ls k}|\{(a,b)\in A\times B\colon  a+\la_jb^{-1}=0\}|
\gs|A|+|B^{-1}|-\deg P.$$
Therefore
$$|C|+\min_{c\in C}|\{(a,b)\in A\times B\colon  ab=c\}|\gs |A|+|B|-\sum_{i=1}^l(m_i+n_i)$$
as desired. \qed

\bigskip
\Ack. The authors are indebted to the referee for his helpful comments.
The revision was done during the second author's visit to
the University of California at Irvine, therefore Sun would like to thank
Prof. Daqing Wan for the invitation.

\enddocument